\documentclass[11pt,reqno]{article}
\usepackage{amsmath,mathrsfs,graphicx,amssymb,amsfonts,color}
\usepackage{cite}

\font\bg=cmbx10 scaled\magstep1

\font\small=cmr8

\newtheorem{newlemma}{{\bf Lemma}}

\newtheorem{newteorem}{{\bf Theorem}}

\newenvironment{teorem}{\begin{newteorem}{\hspace{-0.5
em}{\bf.}}}{\end{newteorem}}

\newtheorem{newkorolari}{{\bf Corollary}}

\newtheorem{newdefine}{{\bf Definition}}

\newtheorem{newquestion}{{\bf Question}}

\newtheorem{newkonjek}{{\bf Conjecture}}

\newtheorem{newexample}{{\bf Example}}

\begin{document}
\tolerance=10000
\baselineskip18truept
\newbox\thebox
\global\setbox\thebox=\vbox to 0.2truecm{\hsize 0.15truecm\noindent\hfill}
\def\boxit#1{\vbox{\hrule\hbox{\vrule\kern0pt
\vbox{\kern0pt#1\kern0pt}\kern0pt\vrule}\hrule}}
\def\qed{\lower0.1cm\hbox{\noindent \boxit{\copy\thebox}}\bigskip}
\def\ss{\smallskip}
\def\ms{\medskip}
\def\bs{\bigskip}
\def\c{\centerline}
\def\nt{\noindent}
\def\ul{\underline}
\def\ol{\overline}
\def\lc{\lceil}
\def\rc{\rceil}
\def\lf{\lfloor}
\def\rf{\rfloor}
\def\ov{\over}
\def\t{\tau}
\def\th{\theta}
\def\k{\kappa}
\def\l{\lambda}
\def\L{\Lambda}
\def\g{\gamma}
\def\d{\delta}
\def\D{\Delta}
\def\e{\epsilon}
\def\lg{\langle}
\def\rg{\rangle}
\def\p{\prime}
\def\sg{\sigma}
\def\ch{\choose}

\newcommand{\ben}{\begin{enumerate}}
\newcommand{\een}{\end{enumerate}}
\newcommand{\bit}{\begin{itemize}}
\newcommand{\eit}{\end{itemize}}
\newcommand{\bea}{\begin{eqnarray*}}
\newcommand{\eea}{\end{eqnarray*}}
\newcommand{\bear}{\begin{eqnarray}}
\newcommand{\eear}{\end{eqnarray}}

\centerline{\Large \bf On the Domination Polynomials of Friendship Graphs}
\vspace{.3cm}

\bigskip

\bs

\baselineskip12truept
\centerline{S. Alikhani$^{a,}${}\footnote{\baselineskip12truept\it\small
Corresponding author. E-mail: alikhani@yazd.ac.ir} J.I. Brown$^b$ and S. Jahari$^a$ }
\baselineskip20truept
\centerline{\it $^a$Department of Mathematics, Yazd University}
\vskip-8truept
\centerline{\it  89195-741, Yazd, Iran}

\centerline{ \it $^{b}$Department of Mathematics and Statistics,
Dalhousie University, Halifax, Canada}

\vskip-0.2truecm
\nt\rule{16cm}{0.1mm}

\nt{\bg ABSTRACT}
\medskip

\baselineskip14truept

\nt{Let $G$ be a simple graph of order $n$.
The {\em domination polynomial} of $G$ is the polynomial
${D(G, x)=\sum_{i=0}^{n} d(G,i) x^{i}}$,
where $d(G,i)$ is the number of dominating sets of $G$ of size $i$.
Let $n$ be any positive integer and  $F_n$ be the {\em Friendship graph} with
$2n + 1$ vertices and $3n$ edges, formed by the join of $K_{1}$ with $nK_{2}$. We study the domination polynomials of this family of graphs, and in particular examine the domination roots of the family, and find the limiting curve for the roots. We also show that for every $n\geq 2$, $F_n$ is not $\mathcal{D}$-unique, that is, there is another non-isomorphic graph with the same domination polynomial.  Also we construct some families of graphs whose real domination roots are only $-2$ and $0$. Finally, we conclude by discussing the domination polynomials of a related family of graphs, the $n$-book graphs $B_n$, formed by joining $n$ copies of the cycle graph $C_4$
with a common edge.}

\ms

\nt{\bf Mathematics Subject Classification:} {\small 05C31, 05C60.}

\nt\rule{16cm}{0.1mm}

\baselineskip20truept

\section{Introduction}

\nt Let $G=(V,E)$ be a simple graph.
For any vertex $v\in V(G)$, the {\it open neighborhood} of $v$ is the
set $N(v)=\{u \in V (G) | uv\in E(G)\}$ and the {\it closed neighborhood} of $v$
is the set $N[v]=N(v)\cup \{v\}$. For a set $S\subseteq V(G)$, the open
neighborhood of $S$ is $N(S)=\bigcup_{v\in S} N(v)$ and the closed neighborhood of $S$
is $N[S]=N(S)\cup S$.
A set $S\subseteq V(G)$ is a {\it dominating set} if $N[S]=V$ or equivalently,
every vertex in $V(G)\backslash S$ is adjacent to at least one vertex in $S$.
The {\it domination number} $\gamma(G)$ is the minimum cardinality of a dominating set in $G$.
For a detailed treatment of domination theory, the reader is referred to~\cite{domination}.

Let ${\cal D}(G,i)$ be the family of dominating sets of a graph $G$ with cardinality $i$ and
let $d(G,i)=|{\cal D}(G,i)|$.
The {\it domination polynomial} $D(G,x)$ of $G$ is defined as
$\displaystyle{D(G,x)=\sum_{ i=\gamma(G)}^{|V(G)|} d(G,i) x^{i}}$ (see \cite{euro,saeid1,kotek}); the polynomial is the generating polynomial for the number of dominating sets of each cardinality.
Similar generating polynomial for other combinatorial sequences, such as independents sets in a graph \cite{bdn,bn,bhn2,brownhickmannowa,chudnovsky,fish90,gut83,gutman,gut90,gut91,gut92b,gut92,gut96}, have attracted recent attention, to name but a few. The algebraic encoding of salient counting sequences allows one to not only develop formulas more easily, but also to prove often unimodality results via the nature of the the roots of the associated polynomials (a well known result of Newton states that if a real polynomial with positive coefficients has all real roots, then the coefficients form a unimodal sequence (see, for example, \cite{comtet}). A root of $D(G, x)$ is called a {\it domination root} of $G$. The set of distinct roots of $D(G, x)$ is denoted by $Z(D(G, x))$ (see \cite{gcom,few,brown}).

Calculating the domination polynomial of a graph $G$ is difficult in general, as the smallest power of a non-zero term is the domination number $\gamma (G)$ of the graph, and determining whether $\gamma (G) \leq k$ is known to be NP-complete \cite{garey}. But for certain classes of graphs, we can find a closed form expression for the domination polynomial. In the next section we will introduce {\em friendship graphs} and calculate their domination polynomials, exploring the nature and location of their roots in conjunction with some outstanding conjectures on domination roots.

\section{Domination polynomials and domination roots of friendship graphs}

\begin{figure}
\begin{center}
\includegraphics[width=6in]{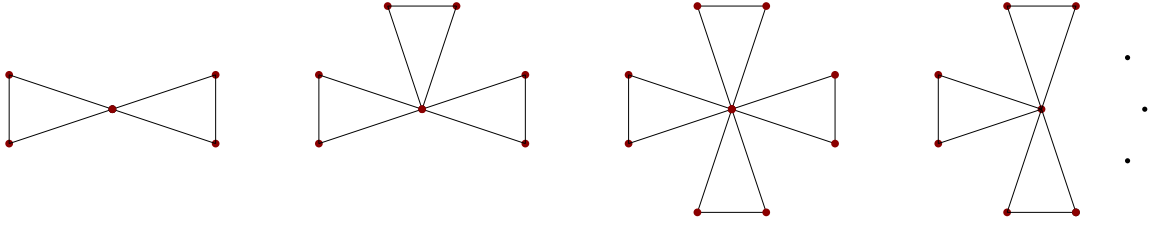}
\caption{Friendship graphs $F_2, F_3, F_4$ and $F_n$, respectively.}
\label{Dutch}
\end{center}
\end{figure}

\nt The friendship (or Dutch-Windmill) graph $F_n$ is a graph that can be constructed by coalescence $n$
copies of the cycle graph $C_3$ of length $3$ with a common vertex. The Friendship Theorem of Paul Erd\"{o}s,
Alfred R\'{e}nyi and Vera T. S\'{o}s \cite{erdos}, states that graphs with the property that every two vertices have
exactly one neighbour in common are exactly the friendship graphs.
Figure \ref{Dutch} shows some examples of friendship graphs.

\begin{teorem} \label{d.p.dutch}
For every $n\in \mathbb{N}$,
$$D(F_n,x) = (2x + x^2)^n +x(1 + x)^{2n}.$$
\end{teorem}
\nt{\bf Proof.} The join $G = G_1 + G_2$ of two graph $G_1$ and $G_2$ with disjoint vertex sets $V_1$ and $V_2$ and
edge sets $E_1$ and $E_2$ is the graph union $G_1\cup G_2$ together with all the edges joining $V_1$ and
$V_2$. An elementary observation is that if $G_1$ and $G_2$ are graphs of orders $n_1$ and $n_2$,
respectively, then
\[ D(G_1 \cup G_2,x) = D(G_1, x) D(G_2, x)\]
and
\[ D(G_1+G_2,x)=\Big((1+x)^{n_1}-1\Big)\Big((1+x)^{n_2}-1\Big)+D(G_1,x)+D(G_2,x). \]
Clearly $D(K_{1},x) = x$ and $D(K_{2},x) = 2x+x^2$, so by  the previous observations,
\begin{eqnarray*}
D(F_n,x) & = & D(k_1+nK_2,x)\\
         & = & (1+x-1)^{1}((1+x)^{2n}-1)+x+(2x+x^{2})^{n}\\
         & = & (2x + x^2)^n +x(1 + x)^{2n}.
\end{eqnarray*}
\hfill \qed

\begin{figure}[h]
\hspace{4cm}
\includegraphics[width=7cm]{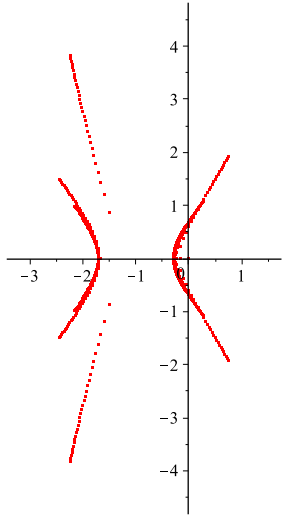}
\caption{\label{figure2'} Domination roots of graphs $F_n$, for $1 \leq n \leq 30$.}
\end{figure}

The domination roots of friendship graphs exhibit a number of interesting properties (see Figure~\ref{figure2'}). Even though we cannot find the roots explicitly, there is much we can say about them.

\subsection{Real domination roots of friendship graphs}

\nt It is known that $-1$ is \underline{not} a domination root as the number of dominating sets in a graph is always odd \cite{brouwer}. On the other hand, of course, $0$ is a domination root of every graph but there are graphs with no nonzero real domination roots. Here we investigate the real domination roots of friendship graphs, and prove first that for any odd natural number $n$, friendship graphs $F_{n}$ have no real domination roots except zero.

\begin{teorem}
For every odd natural number $n$, no nonzero real number is a domination
root of $F_{n}$.
\end{teorem}
\nt{\bf Proof.} By Theorem \ref{d.p.dutch}, for every $n\in
\mathbb{N}$, $D(F_n,x) = (2x + x^2)^n +x(1 + x)^{2n}. $ If
$D(F_n,x)=0$ with $x \neq 0$, then we have
\[
x=-\left( {1-\frac{1}{(1+x)^2}}\right)^n.
\]
We consider three cases, and show in each there is no nonzero solution.
\begin{itemize}
\item  $x > 0:$  Obviously the above equality is true
just for real number $0$, since for nonzero real number the left
side of equality is positive but the right side is negative.
\item $x \leq -2:$ In this case the left side is less than $-2$
and the right side $-\big(1-\frac{1}{(1+x)^2}\big)^n$ is greater
than $-1$, a contradiction.
\item  $-2 < x < 0:$  In this case obviously there are no real solutions $x$ as for odd $n$ and  for every real number $-2<x<0$, the left side of equality is negative but the right side is positive.
\end{itemize}

Thus in any event, there are no nonzero real domination roots of friendship graphs $F_{n}$ where $n$ is odd.
\hfill \qed

We point out that the first two cases also hold when $n$ is even, and hence any real nonzero domination roots of friendship graphs, when $n$ is even, lie in $(-2,0)$, and indeed, it appears that there are always {\em exactly} two real nonzero domination roots in this case. We can show that there are {\em at least} two real nonzero domination roots for $F_{n}$ where $n$ is even: for $n$ even, we see that
\begin{itemize}
\item near but to the left of $0$,
\[ D(F_n,x) = (2x + x^2)^n +x(1 + x)^{2n} = x^{n}(x+2)^n + x(1+x)^{2n} < 0,\]
\item $D(F_n,-1) = (-1)^n  > 0 $, and
\item $D(F_n,-2) = -2(-1)^{2n} < 0.$
\end{itemize}
Hence by the Intermediate Value Theorem, $D(F_n,x)$ has at least two real roots in $(-2,0)$ (with neither being $-1$). Thus the real domination roots of the Friendship graphs are quite different, depending on the parity of $n$.

In fact, for $n\leq 10$, the real roots of $D(F_n,x)$ are (to ten significant digits) shown in Table~\ref{fnrealroots}. The two nonzero real domination roots for $n$ even seem to approach limits, and we will have more to say about this in the next section.

\begin{table}[htp]\label{fnrealroots}
\centering
\fontfamily{ppl}\selectfont
\begin{tabular}{r|l}
$n$ & real domination roots \\  \hline
$1$ & $0$\\
$2$ & $-1.660992532,~-0.1516251043,~0$\\
$3$ & $0$\\
$4$ & $-1.683727169,~-0.2316175850,~0$\\
$5$ & $0$\\
$6$ & $-1.691458147,~-0.2537459684,~0$\\
$7$ & $0$\\
$8$ & $-1.695348455,~-0.2641276712,~0$\\
$9$ & $0$\\
$10$ & $-1.697690028,~-0.2701559954,~0$
\end{tabular}
\caption{Real domination roots of the friendship graph $F_{n}$.}
\label{tab:fnrealroots}
\end{table}

\subsection{Limits of domination roots of friendship graphs}

What about the complex domination roots of friendship graphs? The plot in Figure~\ref{figure2'} suggests that the roots tend to lie on a curve. In order to find the limiting curve, we will need a definition and a well known result.

\begin{newdefine}
If ${f_n(x)}$ is a family of (complex) polynomials, we say that a number $z \in \mathbb{C}$ is a limit of roots of ${f_n(x)}$ if either $f_n(z) = 0$ for all sufficiently large $n$ or z is a limit point of the set $\mathbb{R}({f_n(x)})$, where $\mathbb{R}({f_n(x)})$ is the union of the roots of the $f_n(x)$.
\end{newdefine}

The following restatement of the Beraha-Kahane-Weiss theorem \cite{bkw}  can be found in \cite{brownhickman}.

\begin{teorem}\label{bkw}
Suppose ${f_n(x)}$ is a family of polynomials such that
\begin{eqnarray}
f_n(x) = \alpha_1(x)\lambda_1(x)^n + \alpha_2(x)\lambda_2(x)^n + ... + \alpha_k(x)\lambda_k(x)^n
\end{eqnarray}
where the $\alpha_i(x)$ and the $\lambda_i(x)$ are fixed non-zero polynomials, such that for no pair $i \neq j$ is $\lambda_i(x) \equiv \omega\lambda_j(x)$ for some $\omega \in \mathbb{C}$ of unit modulus. Then $z \in \mathbb{C}$ is a limit of roots of ${f_n(x)}$ if and only if either
\begin{itemize}
\item[(i)] two or more of the $\lambda_i(z)$ are of equal modulus, and strictly greater (in modulus) than the others; or
\item[(ii)] for some $j$, $\lambda_j(z)$ has modulus strictly greater than all the other $\lambda_i(z)$, and $\alpha_j(z) = 0$
\end{itemize}
\end{teorem}

We use Theorem~\ref{bkw} to find the limits of the domination roots of friendship graphs. To do so, we rewrite the domination polynomial
$$D(F_n,x) = (2x + x^2)^n +x(1 + x)^{2n}.$$
of friendship graphs by setting $y = 1+x$. Then we need to consider the limit of roots of
\[f_{n}(y) = (y^{2}-1)^{n}+(y-1)y^{2n},\]
which we rewrite in a form for which we can apply Theorem~\ref{bkw}:
\[f_{n}(y) = (y^{2}-1)^{n}+(y-1)(y^{2})^{n}.\]
We set
\[ \alpha_{1}(y) = 1,~\alpha_{2}(y) = y-1,~\lambda_{1}(y) = y^{2}-1, \mbox{ and } \lambda_{2} = y^{2}.\]
Clearly there is no $\omega \in \mathbb{C}$ of modulus $1$ for which $\lambda_{1} = \omega \lambda_{2}$ (or vice versa), so we can apply  Theorem~\ref{bkw}. Case (ii) is easiest to handle first, as $\alpha_{1}$ is never $0$, and $\alpha_{2} = 0$ if and only if $y = 1$, and in this case $|\lambda_{2}(1)| = |1| > 0 = |\lambda_{1}(1)|$, so we conclude $y =1$ (and hence $x = 0$) is a limit of domination roots of friendship graphs.

The more interesting case is (i), and here we seek all $y$ for which  $|\lambda_{1}(y)| = |\lambda_{2}(y)|$, that is,
\[ |y^{2}-1| = |y^{2}|.\]
To find this curve, let $a = \Re (y)$ and $b = \Im (y)$. Then by substituting in $y = a+ib$ and squaring both sides, we have
\[ (a^{2}-1-b^{2})^{2} + (2ab)^{2} =  (a^{2}-b^{2})^{2} + (2ab)^{2}.\]
This is equivalent to
\[ a^{2} - b^{2} = \frac{1}{2},\]
a hyperbola. Hence, we converting back to variable $x$, we have the following.

\begin{teorem}\label{mainroots}
The limit of domination roots of friendship graphs is $-1$ together with the hyperbola
\[ (\Re x + 1)^{2} + (\Im x)^{2} = \frac{1}{2}.\]
\end{teorem}

\begin{figure}[h]
\hspace{4cm}
\includegraphics[width=7cm]{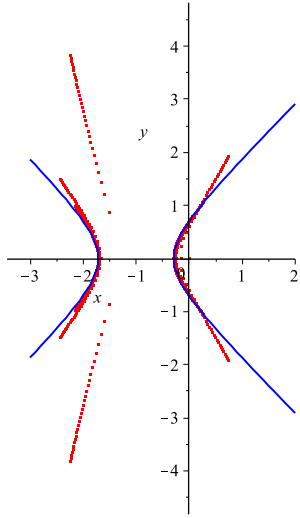}
\caption{\label{hyperbola} Domination roots of graphs $F_n$, for $1 \leq n \leq 30$ along with limiting curve.}
\end{figure}

Figure~\ref{hyperbola} shows the limiting curve. We see that this curve meet the real axis at $\displaystyle{-1-\frac{1}{\sqrt{2}} \approx -1.7071}$ and  $\displaystyle{-1+\frac{1}{\sqrt{2}} \approx -0.2929}$, which agrees well with Table~\ref{tab:fnrealroots}. Also, in \cite{brown} a family of graphs was produced with roots just barely in the right-half plane (showing that not all domination polynomials are stable), but Theorem~\ref{mainroots} provides an explicit family (namely the friendship graphs) whose domination roots have unbounded positive real part.
\ms

\subsection{Uniqueness of domination polynomials of friendship graphs}

\nt Two graphs $G$ and $H$ are said to be {\it dominating equivalent},
or simply ${\cal D}$-equivalent, written $G\sim H$, if
$D(G,x)=D(H,x)$. It is evident that the relation $\sim$ of being
${\cal D}$-equivalence
 is an equivalence relation on the family ${\cal G}$ of graphs, and thus ${\cal G}$ is partitioned into equivalence classes,
called the {\it ${\cal D}$-equivalence classes}. Given $G\in {\cal G}$, let
\[
[G]=\{H\in {\cal G}:H\sim G\}.
\]
We call $[G]$ the equivalence class determined by $G$.
A graph $G$ is said to be dominating unique, or simply
$\mathcal{D}$-unique, if $[G] = \{G\}$, that is, if a graph has the same domination polynomial as $G$, then it must be isomorphic to $G$.

\nt A question of recent interest concerning this equivalence relation $[\cdot]$ asks
which graphs are determined by their domination polynomial. It is known
that cycles \cite{euro} and cubic graphs of order $10$ \cite{cubic} (particularly, the Petersen
graph) are, while if $n\equiv 0 (mod\, 3)$, the paths of order $n$ are not \cite{euro}. In \cite{complete}, authors
completely described the complete $r$-partite graphs which are $\mathcal{D}$-unique. Their results in the
bipartite case, settles in the affirmative a conjecture in \cite{ghodrat}.

\begin{figure}[!h]
\hspace{3cm}
\includegraphics[width=8.5cm,height=2.3cm]{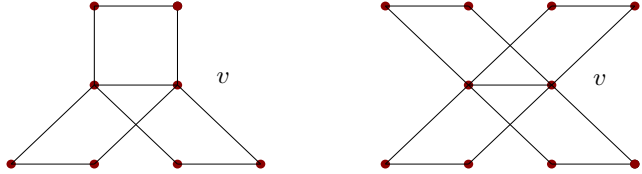}
\caption{ \label{figure6} The book graphs $B_3$ and $B_4$, respectively.}
\end{figure}

What about friendship graphs -- are they $\mathcal{D}$-unique? To answer this question, we introduce a related family of graphs. The $n$-book graph $B_n$ can be constructed by bonding $n$ copies of the cycle graph $C_4$ along a common edge $\{u, v\}$, see Figure \ref{figure6}. We'll now develop a formula for the domination polynomials of book graphs.

%

\nt We begin first with a graph operation.  For two graphs $G = (V,E)$ and $H=(W,F)$, the corona $G\circ H$ is the graph arising from the
disjoint union of $G$ with $| V |$ copies of $H$, by adding edges between
the $i$th vertex of $G$ and all vertices of $i$th copy of $H$ \cite{Fruc}. It is easy to see that the corona operation of two graphs does not have the commutative property.
The following theorem which is for computation of domination
polynomial of corona products of two graphs.

\begin{teorem}\label{theorem7}{\rm \cite{Oper,Kot}}
Let $G = (V,E)$ and $H=(W,F)$ be nonempty graphs of order $n$ and $m$, respectively. Then
\begin{eqnarray*}
D(G\circ H,x) = (x(1 + x)^m + D(H, x))^n.
\end{eqnarray*}
\end{teorem}

\nt The vertex contraction $G/u$ of a graph $G$ by a vertex $u$ is the operation under
which all vertices in $N(u)$ are joined to each other and then $u$ is deleted (see\cite{Wal}).

\nt The following result is  useful for finding the recurrence relations for the  domination polynomials  of  arbitrary graphs.

\begin{teorem}\label{theorem1}{\rm \cite{saeid2,Kot}}
Let $G$ be a graph. For any vertex $u$ in $G$ we have
\[
D(G, x) = xD(G/u, x) + D(G - u, x) + xD(G - N[u], x) - (1 + x)p_u(G, x),
\]
where $p_u(G, x)$ is the polynomial counting the dominating sets of $G - u$ which do not contain any
vertex of $N(u)$ in $G$.
\end{teorem}

\nt Theorem \ref{theorem1} can be used to give a recurrence relation which removes triangles.
Similar to \cite{Kot} we denote  the graph $G\odot u$, graph obtained from $G$ by the removal of all edges between any
pair of neighbors of $u$. Note $u$ is not removed from the graph. The following
recurrence relation is useful on graphs which have many triangles.

\begin{teorem}\label{theorem2}{\rm \cite{Kot}}
Let $G$ be a graph and $u\in V$. Then
\[ D(G, x) = D(G - u, x) + D(G\odot u, x) - D(G \odot u - u, x). \]
\end{teorem}

%

\nt We are now ready to give a formula for the domination polynomial of $B_n$.

\begin{teorem}\label{theorem11}
\nt For every $n \in \mathbb{N}$,
\[ D(B_n,x)=(x^2+2 x)^n(2x+1) + x^2(x+1)^{2n}- 2x^n.\]
\end{teorem}
\nt{\bf Proof.}
\nt Consider graph  $B_{n}$ and a vertex v in the common edge (see Figure 5). By Theorems ~\ref{theorem1} we have:
\begin{eqnarray*}
D(B_n, x)&=& x D(B_n/v, x) + D(B_n - v, x) + x D(B_n - N[v], x) - (1 + x)p_v(B_n, x)\\
&=& x D(B_n/v, x) + D(B_n-v, x) + x(D(nK_1,x)) - (1 + x)x^n\\
&=&x D(B_n/v, x) + D(B_n-v, x) - x^n.
\end{eqnarray*}

\begin{figure}[!h]
\hspace{2.3cm}
\includegraphics[width=5cm,height=2.9cm]{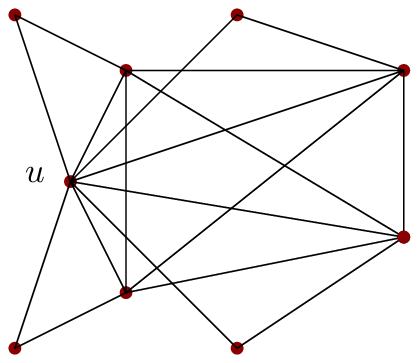}
\hspace{1.5cm}
\includegraphics[width=5cm,height=2.9cm]{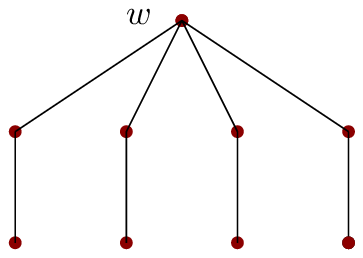}
\caption{ \label{figure8} Graphs $B_4/v$ ~and ~ $B_4 - v$,~ respectively.}
\end{figure}

\begin{figure}[!h]
\hspace{2.9cm}
\includegraphics[width=4.5cm,height=2.7cm]{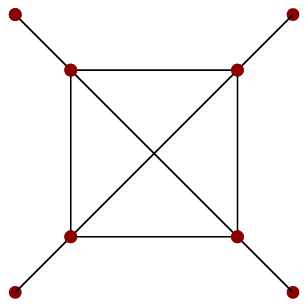}
\hspace{1.5cm}
\includegraphics[width=3.5cm,height=2.3cm]{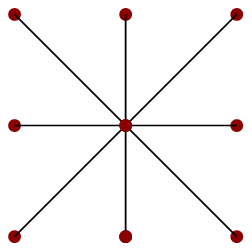}
\caption{ \label{figure9} Graphs $B_4/v-u$ ~and ~ $B_4 - N[v]\odot u$,~ respectively.}
\end{figure}

\nt Now  we use Theorems~\ref{theorem2} to obtain the domination polynomial of the graph $B_n/v$
 (see Figure~\ref{figure8}). We have
\[D(B_n/v, x) = D((B_n/v) - u, x) + D((B_n/v) \odot u, x) - D(2nK_1,x),\]
\nt where $(B_n/v) - u=K_n\circ K_1$ and $(B_n/v) \odot u=K_{1,2n}$ (see Figure~\ref{figure9}).

\nt Using Theorems~\ref{theorem2}, we deduce that,
$ D(B_n/v,x) =(2x+x^2)^{n} +x(x+1)^{2n}.$
 Also we use Theorems ~\ref{theorem7} and ~\ref{theorem1}  to obtain the domination polynomial of the graph $B_n-v$ (see Figure~\ref{figure8}). Hence $D(B_n-v, x) = xD((B_n-v)/w, x) + D(K_2, x)^n - x^n,$
\nt where $(B_n-v)/w=K_n\circ K_1$. So $D(B_n-v, x) = (2x+x^2)^n(x+1) - x^n$.
  Note that in this case $p_v(B_n, x)=p_w(B_n-v, x)=x^n.$  Consequently,
\begin{eqnarray*}
D(B_n, x)&=&x ((2x+x^2)^{n} +x(x+1)^{2n}) + (2x+x^2)^n(x+1) - x^n - x^n\\
&=& (x^2+2 x)^n(2x+1) + x^2(x+1)^{2n}- 2x^n.\quad\qed
\end{eqnarray*}

\begin{teorem}
For each natural number $n\geq 2,$ the friendship  graph $F_n$ is not $\mathcal{D}$-unique, as  $F_n$ and $B_n/v$ have the same domination polynomial.
\end{teorem}
\nt{\bf Proof.}  In the proof of Theorem \ref{theorem11}, we proved that $ D(B_n/v,x) =(2x+x^2)^{n} +x(x+1)^{2n}.$ Therefore
 $D(F_n,x)=D(B_n/v,x)$. Since $F_n$ is not isomorphic to $B_n/v$, for each natural number $n\geq 2$, so the friendship graphs are not $\mathcal{D}$-unique and  $[F_n]\supseteq \{F_n,B_n/v\}$.\quad\qed

\section{Open Problems}

The results of the previous section show that even if we can find an explicit formula for the domination polynomial of a graph, there are still interesting, difficult problems concerning the roots.  With regards to friendship graphs, we pose the following:

\begin{newquestion}
For $n$ even, does $F_{n}$ have exactly three real roots?
\end{newquestion}

\begin{newquestion}
What is a good upper bound on the modulus of the roots of  $F_{n}$?
\end{newquestion}

\begin{figure}[h]\label{bookroots}
\hspace{4cm}
\includegraphics[width=7cm]{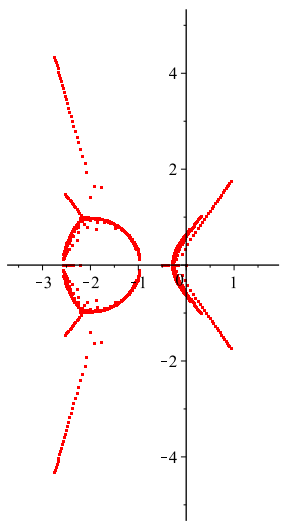}
\caption{\label{hyperbola} Domination roots of graphs $B_n$, for $1 \leq n \leq 30$.}
\end{figure}

Some calculations seem to indicate that the moduli of the roots, while going off to infinity (by Theorem~\ref{mainroots}), do so quite slowly, perhaps like $\ln n$.
The book graphs indeed have a more interesting formula than friendship graphs. Figure~\ref{bookroots} shows the domination roots of $n$-book graph for $n \leq 30$. Questions about the real roots, the limit of the roots, bounding the moduli of the roots can be asked as well. (We remark that using Theorem~\ref{bkw}, we can show that the limit of the roots is the circle $|x+2|= 1$ with real part at least $\displaystyle{-\frac{3}{2}-\frac{\sqrt{2}}{2}}$, the portion of the hyperbola $ (\Re x + 1)^{2} + (\Im x)^{2} = \frac{1}{2}$ in the right half-plane, plus the portion of the curve $|x+1|^{2} = |x|$ with real part at most  $\displaystyle{-\frac{3}{2}-\frac{\sqrt{2}}{2}}$.)

\begin{newquestion}
What can be said about the domination roots of book graphs?
\end{newquestion}

\nt Along these lines, there is a conjecture  which states that, the set of integer domination roots of any graphs
is a subset of $\{-2,0\}$ (\cite{few}). Now we show that there are infinite families of graphs, based on friendship and book graphs, whose their domination polynomials
have real roots $-2$ and $0$.

\begin{teorem}
\begin{enumerate}
\item [(i)] For every  odd natural number $n$, the only nonzero real domination  root of  $B_n\circ F_n$  is $-2$.

\item [(ii)]  For every  even natural number $n$, the only nonzero real domination root of $B_n\circ F_{n+1}$ and $B_{n+1}\circ F_n$  is $-2$.
\end{enumerate}
\end{teorem}

\nt{\bf Proof.} (i) By theorems \ref{d.p.dutch} and \ref{theorem7}  for any odd natural number $n$ and arbitrary  graph $G$ we have,
\begin{eqnarray*}\label{eq2}
D(G\circ F_n,x)&=&\Big(x(1+x)^{2n+1}+(2x+x^2)^n+x(1+x)^{2n}\Big)^{|V(G)|}\nonumber\\
 &=&\Big(x[(1+x)^{2n+1}+x^{n-1}(2+x)^n+(1+x)^{2n}]\Big)^{|V(G)|}\nonumber\\
  &=&\Big(x[(1+x)^{2n}(1+x+1)+x^{n-1}(2+x)^n]\Big)^{|V(G)|}\nonumber\\
    &=&\Big(x(2+x)[(1+x)^{2n}+x^{n-1}(2+x)^{n-1}]\Big)^{|V(G)|}\nonumber\\
        &=&(x(2+x))^{|V(G)|}\Big([(1+x)^{2n}+(2x+x^2)^{n-1}]\Big)^{|V(G)|}.
\end{eqnarray*}
        Now we prove that, for each  odd natural $n$, $f_n(x)=(2x + x^2)^{n-1} +(1 + x)^{2n}$ have no   real roots.  If
$f_n(x)=0$, then we have
\[((1+x)^2-1)^{n-1}=-(1+x)^{2n}. \]
Obviously the above equality is not true for any real number. Because for odd $n$ and  for every real number, the left side of equality is positive but the right side is negative.

\nt (ii) Proof is similar to proof of Part (i).\quad\qed

Along the same lines, we can  show:

\begin{teorem}
\begin{enumerate}
\item[(i)] Every graph $H$ in
the family $\{G\circ K_{2n}, (G\circ K_{2n})\circ K_{2n}, ((G\circ K_{2n})\circ K_{2n})\circ K_{2n},\cdots \}$ does not have real domination roots, except zero.

\item[(ii)]
Every graph $H$ in
the family $\{G\circ K_{2n+1}, (G\circ K_{2n+1})\circ K_{2n+1}, ((G\circ K_{2n+1})\circ K_{2n+1})\circ K_{2n+1},\cdots \}$ does not have real domination roots, except $\{-2,0\}$.

\item[(iii)] Every graph $H$ in
the family $\{G\circ B_2, (G\circ B_2)\circ B_2, ((G\circ B_2)\circ B_2)\circ B_2,\cdots \}$ does not have real domination roots, except zero.
\end{enumerate}
\end{teorem}

We end off on a final problem.

\begin{newquestion}
Is $-2$ the only possible nonzero integer domination root?
\end{newquestion}


\end{document}